\documentclass[reqno,centertags,a4paper]{amsart}
\usepackage{amssymb,amsmath,amsthm}
\usepackage{hyperref}
\usepackage{times}

\newcommand{\R}{\mathbb{R}}

\theoremstyle{plain}
\newtheorem{theorem}{Theorem}[section]
\newtheorem{lemma}{Lemma}[section]

\theoremstyle{definition}
\newtheorem{definition}{Definition}[section]
\theoremstyle{remark}
\newtheorem{remark}{Remark}

\newcommand{\eps}{\varepsilon}
\newcommand{\lb}{\langle}
\newcommand{\rb}{\rangle}

\author[S. Herr]{Sebastian Herr}

\title{An improved bilinear estimate for Benjamin-Ono type equations}


\address{Fachbereich Mathematik, Universit\"at Dortmund, 44221 Dortmund, Germany.}
\email{sebastian.herr@math.uni-dortmund.de}

\begin{document}

\begin{abstract}
A bilinear estimate in Fourier restriction norm spaces with applications to the Cauchy problem
\begin{align*}
u_t - |D|^{\alpha}u_x + uu_x &=0 \quad \text{in }(-T,T) \times \R\\
u(0)&=u_0
\end{align*}
is proved, for $1< \alpha <2$.
As a consequence, local  well-posedness in $H^s(\R) \cap \dot{H}^{-\omega}(\R)$ follows
for $$s >-\frac{3}{4}(\alpha-1) \,\text{ and }\, \omega=1/\alpha-1/2$$
This extends to global well-posedness for all $s \geq 0$.
\end{abstract}


\maketitle

\section{Introduction}\label{sect:intro}
We consider the Cauchy problem
\begin{equation}\label{eq:bo_alpha}
\begin{split}
u_t - |D|^{\alpha}u_x + uu_x &=0 \quad \text{in }(-T,T) \times \R\\
u(0)&=u_0
\end{split}
\end{equation}
for $1< \alpha <2$ and we are interested in well-posedness results in low regularity
Sobolev spaces.

Our aim is to give an improvement to our previous results \cite{H05a},
where we proved that the Cauchy problem is locally well-posed in
$H^s(\R)\cap \dot{H}^{-\omega}(\R)$ for $s \geq 1-\alpha/2$
and $\omega=\frac{1}{\alpha}-\frac{1}{2}$. Due to the conserved Hamiltonian
this also implied global well-posedness for $s\geq \frac{\alpha}{2}$.
Moreover, using the counterexamples found by Molinet, Saut and
Tzvetkov in \cite{MST01} it was shown that the condition on the low frequencies
is sharp in the sense that for $\omega<\frac{1}{\alpha}-\frac{1}{2}$ the flow
map fails to be $C^2$. For further references and results we refer the reader
to the works of Colliander, Kenig and Staffiliani \cite{CKS02} and Kenig and
Koenig \cite{KK03} and the introduction of \cite{H05a}.
After \cite{H05a} was completed, we learned that these
results were also an improvement to a similar approach by Molinet and Ribaud \cite{MR03b}.

Here, by a refined bilinear estimate we observe that 
the same local well-posedness result holds for all $s >-\frac{3}{4}(\alpha-1)$, which immediately
implies global well-posedness for all $s\geq 0$.

Our analysis includes the range $1<\alpha<2$, without the endpoints
$\alpha=1,2$. Very recently, Kenig and Ionescu \cite{IK05} studied global
well-posedness of the Benjamin-Ono equation ($\alpha=1$) for real valued data in $L^2$ (see
also \cite{BP05,T03}). 
We observe that in the limit for $\alpha \to 2$ our lower bound on $s$ tends to
$-\frac{3}{4}$ which coincides with the results of Kenig, Ponce and Vega \cite{KPV96}
for the Korteweg-de Vries equation and the low frequency condition disappears. In the
limit for $\alpha \to 1$ the lower bound for $s$ tends to $0$ and $\omega \to \frac{1}{2}$.
We believe that the lower bound for $s$ is optimal, but this is work in
progress.

\subsection*{Acknowledgments}\label{subsect:ack}
The author is grateful to M. Hadac and H. Koch for discussions on the
subject. Moreover, the author would like to thank A. Ionescu and C.E. Kenig for
interesting remarks on the Benjamin-Ono case.

\section{Notation and Definition of the Spaces}\label{sect:not}
Let $\mathcal{S}(\R^n)$ be the space of Schwartz
functions on $\R^n$ and define the Fourier transform by
$$\mathcal{F}f(\xi)=\widehat{f}(\xi)=(2\pi)^{-\frac{n}{2}}\int_{\R^n} e^{-i x
  \cdot \xi} f(x) \, dx$$
The partial Fourier transform w.r.t. $t \in \R$ ($x \in \R$)
will be denoted by $\mathcal{F}_t$ ($\mathcal{F}_x$).
$|D|^{s}$ denotes the Fourier multiplier operator with
$\mathcal{F}|D|^{s}v(\xi)=|\xi|^{s}\mathcal{F}v(\xi)$, and $J^s$ is the operator
with symbol $\lb \xi \rb^{s}$.

We write
$$W_\alpha(t): H^{(s,\omega)} \to H^{(s,\omega)}, \mathcal{F}_x W_\alpha(t)u_0 (\xi)
=e^{it\xi|\xi|^{\alpha}} \mathcal{F}_x u_0(\xi)$$
for the solution operator of the linear homogeneous problem, which defines a
unitary group on $H^{(s,\omega)}$.

Throughout this work let $\psi \in C^{\infty}_{0}([-2,2])$ be a nonnegative,
symmetric function with $\psi|_{[-1,1]}\equiv 1$ and let
$\psi_T(t):=\psi(t/T)$.

We use the same spaces as in \cite{H05a}
\begin{definition}\label{def:init_space}For $s \geq 0$ and
$0\leq \omega<\frac{1}{2}$ we define the Sobolev space
$H^{(s,\omega)}$ as the completion of $\mathcal{S}(\R)$
with respect to the norm
\begin{equation}\label{eq:init_space_norm}
\|u\|^2_{H^{(s,\omega)}}:= \int_{\R} \lb \xi \rb
^{2s+2\omega}|\xi|^{-2\omega}|\widehat{u}(\xi)|^2 \, d\xi.
\end{equation}
\end{definition}
Sometimes it is convenient to identify $H^{(s,\omega)}$ and $H^s(\R) \cap \dot{H}^{-\omega}(\R)$.
Our resolution space, a variant of the Bourgain spaces introduced in \cite{Bo93}, will be
\begin{definition}\label{def:res_space} For $0 \leq \omega <\frac{1}{2}$ and $s,b \in \R$ we define
the space $X_{s,\omega,b}$ as the completion of $\mathcal{S}(\R^2)$ with respect
to the norm
\begin{equation}\label{eq:res_space_norm}
\|u\|^2_{X_{s,\omega,b}}:=\int_{\R^2}
|\xi|^{-2\omega}\lb\xi\rb^{2s-2\alpha\omega} \lb |\tau|+|\xi|^{1+\alpha}
\rb^{2\omega}\lb\tau-\xi|\xi|^{\alpha}\rb^{2b}|\mathcal{F}u(\tau,\xi)|^2\,
d\tau d\xi. 
\end{equation}
For $T>0$ we define the restriction norm space 
$$X^T_{s,\omega,b}:=\{u|_{[-T,T]} \mid u \in X_{s,\omega,b} \}$$
with norm
$$\|u\|_{X^T_{s,\omega,b}}=\inf\{\|\widetilde{u}\|_{X_{s,\omega,b}} \mid u=
\widetilde{u}|_{[-T,T]}, \, \widetilde{u}\in X_{s,\omega,b}  \}.$$
\end{definition}
\section{Main Results}\label{sect:main}
Our aim is to prove the following bilinear estimate.
\begin{theorem}\label{thm:bil_est}
Let $1< \alpha <2$, $s\geq s_0 >-\frac{3}{4}(\alpha-1)$ and $\omega=\frac{1}{\alpha}-\frac{1}{2}$.
There exists $b' >-\tfrac{1}{2}$ and $b \in
(\tfrac{1}{2},b'+1)$ such that
\begin{equation}\label{eq:bil_est}
\|\partial_x(u_1u_2)\|_{X_{s,\omega,b'}}
\leq c \|u_1\|_{X_{s,\omega,b}}\|u_2\|_{X_{s_0,\omega,b}}
+\|u_1\|_{X_{s_0,\omega,b}}\|u_2\|_{X_{s,\omega,b}}
\end{equation}
for all $u_1,u_2 \in \mathcal{S}(\R^2)$.
\end{theorem}
This leads to local well-posedness by an application of the contraction
mapping principle in a straightforward way.
For the general outline of the proof we refer the reader to e.g. \cite{Bo93,G96,KPV96}. The
minor modifications of these arguments in the $X_{s,\omega,b}$ spaces
are carried out in detail in our previous work \cite{H05a}.
\begin{theorem}\label{thm:main_local} Let $1<\alpha<2$ and
$\omega=\frac{1}{\alpha}-\frac{1}{2}$. Then,
for $s \geq s_0 >-\frac{3}{4}(\alpha-1)$ there exists $b>\frac{1}{2}$ and a non-increasing function $T:
(0,\infty) \to (0,\infty)$, such that for
any $u_0 \in H^{(s,\omega)}$ and  $T=T(\|u_0\|_{H^{(s_0,\omega)}})$,
there exists a solution
$$u \in X^T_{s,\omega,b} \subset C\big([-T,T],H^{(s,\omega)}\big)$$
of the Cauchy problem
\begin{equation*}
\begin{split}
u_t - |D|^{\alpha}u_x + uu_x &=0 \quad \text{in }(-T,T) \times \R\\
u(0)&=u_0
\end{split}
\end{equation*}
which is unique in the class
of $X^T_{s_0,\omega,b}$ solutions.
Moreover, for any $r>0$ there exists $T=T(r)$, such that for
$B=\{v_0 \in H^{(s,\omega)}\mid
\|v_0\|_{H^{(s_0,\omega)}} \leq r\}$
the flow map 
$$
F: H^{(s,\omega)} \supset B  \to C\big([-T,T],H^{(s,\omega)}\big) \cap
X^{T}_{s,\omega,b} \quad,\; u_0 \mapsto u
$$
is analytic.
\end{theorem}
\begin{remark}\label{rem:solution}
Here, solution always means fixed point of (an extension of) the operator
$$\Phi_T(u)(t)=\psi(t)W_\alpha(t)u_0-\frac{1}{2}\psi_T(t)\int_0^tW_\alpha(t-t')\partial_x(u^2)(t')\,dt' $$
in $X_{s,\omega,b}$. These solutions are
solutions in the sense of distributions at least\footnote{Even for $s<0$
one can still use some smoothing properties to verify this} for $s\geq 0$.
\end{remark}
Together with the a priori bound from Lemma \ref{lem:apriori_bound} this also shows the following
\begin{theorem}\label{thm:main_global} Let $1<\alpha<2$ and
$\omega=\frac{1}{\alpha}-\frac{1}{2}$ as well as $s_0>-\frac{3}{4}(\alpha-1)$. Then,
for $s \geq 0$ there exists $b>\frac{1}{2}$, such that for every $T>0$
and real valued $u_0 \in H^{(s,\omega)}$
there exists a real valued solution
$$u \in X^T_{s,\omega,b} \subset C\big([-T,T],H^{(s,\omega)}\big)$$
of the Cauchy problem
\begin{equation*}
\begin{split}
u_t - |D|^{\alpha}u_x + uu_x &=0 \quad \text{in }(-T,T) \times \R\\
u(0)&=u_0
\end{split}
\end{equation*}
which is unique in $X^T_{s_0,\omega,b}$.
Moreover, the flow map 
$$
F: H^{(s,\omega)} \to C\big([-T,T],H^{(s,\omega)}\big) \cap
X^{T}_{s,\omega,b} \quad,\; u_0 \mapsto u
$$
is real analytic.
\end{theorem}
\section{Preparatory Lemmata}\label{sect:prep}
In this section we will summarize our main tools for the proof of the bilinear estimate.
First, we recall the $L^4_tL^{\infty}_x$ Strichartz estimate.
\begin{lemma}\label{lem:strichartz}
For $b>\frac{1}{2}$ we have
\begin{equation}
\|J^{\frac{\alpha-1}{4}}u\|_{L^4_tL^{\infty}_x} \leq c \|u\|_{X_{0,0,b}}
\label{eq:str_with_deriv}
\end{equation}
\begin{proof}
From \cite{KPV91b} Theorem 2.1, we know that
\begin{equation*}
\||D|^{\frac{\alpha-1}{4}}W_{\alpha}(t)u_0\|_{L^4_tL^{\infty}_x} \leq c
\|u_0\|_{L^2}
\end{equation*}
By the general properties of Bourgain spaces, see e.g. \cite{G96} Lemme 3.3, the
estimate
\begin{equation}\label{eq:kpv_xsb}
\||D|^{\frac{\alpha-1}{4}}u\|_{L^4_tL^{\infty}_x} \leq c \|u\|_{X_{0,0,b}}
\end{equation}
follows.
By smooth cutoffs in frequency, we split $u$ into a low frequency part $u^{low}$
with
$$\mathcal{F}u^{low}(\tau,\xi)=\psi(\xi)\mathcal{F}u(\tau,\xi)$$
and a high frequency part $u^{high}:=u-u^{low}$.
Then,
$$
\|J^{\frac{\alpha-1}{4}}u\|_{L^4_tL^{\infty}_x} \leq 
\|J^{\frac{\alpha-1}{4}}u^{low}\|_{L^4_tL^{\infty}_x}+
\|J^{\frac{\alpha-1}{4}}u^{high}\|_{L^4_tL^{\infty}_x}
$$
By an application of the Sobolev inequality, the  first part is bounded by
$$c\|J^{\frac{\alpha+1}{4}+\eps}
u^{low}\|_{L^4_tL^{2}_x}
\leq c
\|u\|_{L^4_tL^{2}_x}\leq c
\|u\|_{X_{0,0,b}}$$
whereas the second term is
bounded by
$$
c\||D|^{-\frac{\alpha-1}{4}}J^{\frac{\alpha-1}{4}}
u^{high}\|_{X_{0,0,b}} \leq c\|u\|_{X_{0,0,b}}
$$
due to \eqref{eq:kpv_xsb}, which gives the desired estimate.
\end{proof}
\end{lemma}
The next Lemma contains a bilinear Strichartz type estimate in the spirit of \cite{G01b,G01a}.
For the proof we refer to our previous work \cite{H05a}.
\begin{lemma}\label{lem:bil_str}
We define the bilinear operator $I_{\ast}^{s}$ via
\begin{align*}
\mathcal{F} \, I_{\ast}^{s}(u_1,u_2)(\tau,\xi)
=\int_{\genfrac{}{}{0pt}{}{\xi=\xi_1+\xi_2}{\tau=\tau_1+\tau_2}} 
\left||\xi_1|^{2s}-|\xi_2|^{2s}\right|^{\frac{1}{2}}
\mathcal{F}u_1(\tau_1,\xi_1)\mathcal{F}u_2(\tau_2,\xi_2)
\,d\tau_1 d\xi_1
\end{align*}
For $b>\frac{1}{2}$
\begin{equation}\label{eq:bil_str}
\left\|I_{\ast}^{\frac{\alpha}{2}}(u_1,u_2)\right\|_{L^2_{xt}}
\leq c \|u_1\|_{X_{0,0,b}}\|u_2\|_{X_{0,0,b}} \quad, u_1,u_2 \in X_{0,0,b}
\end{equation}
Moreover, we define $K_{\ast}^{\frac{\alpha}{2}}$ as
\begin{equation*}
\mathcal{F} \,
K_{\ast}^{\frac{\alpha}{2}}(u_1,u_2)(\tau,\xi)=\int_{\genfrac{}{}{0pt}{}{\xi=\xi_1+\xi_2}{\tau=\tau_1+\tau_2}} 
\left||\xi|^{\alpha}-|\xi_1|^{\alpha}\right|^{\frac{1}{2}}
\mathcal{F}\overline{u}_1(\tau_1,\xi_1)\mathcal{F}u_2(\tau_2,\xi_2)
\,d\tau_1 d\xi_1
\end{equation*}
$K_{\ast}^{\frac{\alpha}{2}}$ 
is the formal adjoint of $u_2 \mapsto I^{\frac{\alpha}{2}}_{\ast}(u_1,u_2)$
with respect to $L^2_{xt}$ and 
for $b >\frac{1}{2}$
\begin{equation}\label{eq:dual_bil_str}
\left\|K_{\ast}^{\frac{\alpha}{2}}(u_1,u_2)\right\|_{X_{0,0,-b}}
\leq c \|u_1\|_{X_{0,0,b}}\|u_2\|_{L^2_{xt}} \quad, u_1 \in X_{0,0,b}\;, u_2 \in L^2_{xt}
\end{equation}
\end{lemma}
Finally, we note the elementary resonance relation, which is crucial to exploit the weights in our resolution
space.
\begin{lemma}\label{lem:reson}
Let $1 < \alpha <2$. Define
\begin{equation*}
h(\xi_1,\xi_2,\xi)=\xi|\xi|^{\alpha}-\xi_1|\xi_1|^{\alpha}-\xi_2|\xi_2|^{\alpha}
\end{equation*}
Then, for $\xi=\xi_1+\xi_2$ it holds that
\begin{equation}\label{eq:res_est}
|h(\xi_1,\xi_2,\xi)| \geq c |\xi_{\min}||\xi_{\max}|^{\alpha},
\end{equation}
with $|\xi_{\min}|:=\min\{|\xi_1|,|\xi_2|,|\xi|\}$
and $|\xi_{\max}|:=\max\{|\xi_1|,|\xi_2|,|\xi|\}$.
\end{lemma}

\section{Proof of the bilinear estimate}\label{sect:proof_bilinear}
Let us fix notation.
We define $\sigma=|\tau|+|\xi|^{1+\alpha}$
and $\sigma_i=|\tau_i|+|\xi_i|^{1+\alpha}$ as well as
$\lambda=\tau-\xi|\xi|^{\alpha}$ and
$\lambda_i=\tau_i-\xi_i|\xi_i|^{\alpha}$.
Moreover, we set
$$f_i(\tau_i,\xi_i)=|\xi_i|^{-\omega}\lb\xi_i\rb^{s-\alpha\omega}
\lb\lambda_i\rb^{b}\lb\sigma_i\rb^{\omega}\mathcal{F}u_i(\tau_i,\xi_i)$$
and
$$\mathcal{F}v_i(\tau_i,\xi_i):=f_i(\tau_i,\xi_i)\lb\lambda_i\rb^{-b}.$$
We use the notation
$$\int_{\ast}g(\tau_1,\xi_1)h(\tau_2,\xi_2)
:=\int_{\genfrac{}{}{0pt}{}{\xi=\xi_1+\xi_2}{\tau=\tau_1+\tau_2}}
g(\tau_1,\xi_1)h(\tau_2,\xi_2) \,d\tau_1d\xi_1$$

We first consider the case $s=s_0=-\frac{3}{4}(\alpha-1)+\eps$ for small $\eps>0$.
Our goal is to bound
\begin{equation*}
\|\partial_x(u_1u_2)\|_{X_{s,\omega,b'}}=\left\|
|\xi|^{1-\omega}\lb\xi\rb^{s-\alpha\omega}\lb\lambda\rb^{b'}\lb\sigma\rb^{\omega}
\int_\ast\prod\limits_{i=1}^2
\frac{|\xi_i|^{\omega}\lb\xi_i\rb^{\alpha\omega-s}f_i(\tau_i,\xi_i)}
{\lb\lambda_i\rb^{b}\lb\sigma_i\rb^{\omega}}
\right\|_{L^2_{\tau,\xi}}
\end{equation*}
by the product of the $L^2$ norms of the $f_i$, where we may assume that $0 \leq f_i \in
\mathcal{S}(\R^2)$.

Due to the symmetry in $\xi_1,\xi_2$ it
suffices to consider the subregion of the domain of integration where
$|\xi_1|\leq |\xi_2|$. By the convolution constraint $\xi=\xi_1+\xi_2$
we then have $|\xi|\leq 2|\xi_2|$. This region is splitted again into
\begin{enumerate}
\item Region $D_1$: $4|\xi_1| \leq |\xi_2|$.
  There, $|\xi_1|\leq \tfrac{1}{4} |\xi_2| \leq \tfrac{1}{3} |\xi| \leq
  \tfrac{2}{3} |\xi_2|$.
\item Region $D_2$: $|\xi_1|\leq |\xi_2|\leq 4|\xi_1|$. There, $|\xi| \leq 2|\xi_2|$, $|\xi|\leq 5|\xi_1|$.
\end{enumerate}
Let $A,A_1,A_2$ be subregions of the domain of integration, such that in
$A$ we have $\lb\lambda\rb \geq  \lb\lambda_1\rb ,\lb\lambda_2\rb$, in $A_1$
we have
$\lb\lambda_1\rb \geq  \lb\lambda\rb ,\lb\lambda_2\rb$ and in
$A_2$ the inequalities
$\lb\lambda_2\rb \geq  \lb\lambda\rb ,\lb\lambda_1\rb$ hold.

We first consider the region $D_1$ and subdivide it into two parts
$D_1=D_{11}\cup D_{12}$, where in $D_{11}$ we have $|\xi_1|\leq 2$ and in $D_{12}$
we have $|\xi_1|\geq 2$.
In $D_1$ we see by Lemma \ref{lem:reson}
$$
|\lambda-\lambda_1-\lambda_2|
=
|h(\xi_1,\xi_2,\xi)|\geq c |\xi_1||\xi|^\alpha
$$
because $|\xi_1|=|\xi_{\min}|$ and $|\xi|\leq 2|\xi_{\max}|$.

Now we start the analysis in the subregion $D_{11}$ where the arguments
remain close to those in \cite{H05a}.
We exploit
$$|\xi|^{1-\frac{\alpha}{2}} =|\xi|^{\alpha\omega}\leq c |\xi_1|^{-\omega}
(\chi_A\lb\lambda\rb^{\omega}+\chi_{A_1}\lb\lambda_1\rb^{\omega}+\chi_{A_2}\lb\lambda_2\rb^{\omega}).$$
Therefore in $D_{11}$ the bilinear estimate follows from
\begin{equation}\label{eq:case11}
\sum_{k=0}^2 \|J_{11,k}\|_{L^2} \leq c \prod\limits_{i=1}^2 \|f_i\|_{L^2},
\end{equation}
where
$$
J_{11,0}=\int_{\ast}\chi_{D_{11}\cap A}
|\xi|^{\frac{\alpha}{2}-\omega}\lb\xi\rb^{s-\alpha\omega}\lb\lambda\rb^{b'+\omega}
\lb\sigma\rb^{\omega}|\xi_2|^{\omega}\prod\limits_{i=1}^2
\frac{f_i(\tau_i,\xi_i)\lb\xi_i\rb^{\alpha\omega-s}}
{\lb\lambda_i\rb^{b}\lb\sigma_i\rb^{\omega}}
$$
and for $k=1,2$
$$
J_{11,k}=\int_{\ast}\chi_{D_{11}\cap A_k}
|\xi|^{\frac{\alpha}{2}-\omega}\lb\xi\rb^{s-\alpha\omega}\lb\lambda\rb^{b'}\lb\sigma\rb^{\omega}
|\xi_2|^{\omega}\lb\lambda_k\rb^{\omega}\prod\limits_{i=1}^2
\frac{f_i(\tau_i,\xi_i)\lb\xi_i\rb^{\alpha\omega-s}}
{\lb\lambda_i\rb^{b}\lb\sigma_i\rb^{\omega}}
$$
We observe that in $D_{11}$
\begin{equation}
  \label{eq:freq_bound}
\lb\xi_2\rb^{\alpha\omega-s}\lb\xi\rb^{s-\alpha\omega} \leq
c \text{ and } \lb\xi_1\rb^{\alpha\omega-s}\leq c  
\end{equation}
In addition, we use $b'+\omega \leq 0$ and $|\xi_2|^{\omega} \leq
c|\xi|^{\omega}$ to show that
$$
\|J_{11,0}\|_{L^2} \leq c
\left\|\int_{\ast}\chi_{D_{11}\cap A}
|\xi|^{\frac{\alpha}{2}}
\lb\sigma\rb^{\omega}\prod\limits_{i=1}^2
\frac{f_i(\tau_i,\xi_i)}
{\lb\lambda_i\rb^{b}\lb\sigma_i\rb^{\omega}}\right\|_{L^2}
$$
Because of the convolution constraint $(\tau,\xi)=(\tau_1,\xi_1)+(\tau_2,\xi_2)$
we also have
\begin{equation}\label{eq:sigma}
\frac{\lb\sigma\rb}{\lb\sigma_1\rb\lb\sigma_2\rb}
\leq c\frac{1}{\min_{i=1,2}\lb\sigma_i\rb}\leq c
\end{equation}
which implies
$$
\|J_{11,0}\|_{L^2}\leq c\left\|
\int_{\ast}\chi_{D_{11}\cap A}
|\xi|^{\frac{\alpha}{2}}
\prod\limits_{i=1}^2
\frac{f_i(\tau_i,\xi_i)}
{\lb\lambda_i\rb^{b}}
\right\|_{L^2}
$$
We observe that in $D_{11}$
$$
|\xi|^{\frac{\alpha}{2}}\leq c |
|\xi_2|^{\alpha}-|\xi_1|^{\alpha}|^{\frac{1}{2}},
$$
such that with
\eqref{eq:bil_str}
\begin{align*}
\|J_{11,0}\|_{L^2}&\leq c \left\|
\int_{\ast}
| |\xi_2|^{\alpha}-|\xi_1|^{\alpha}|^{\frac{1}{2}}
\prod\limits_{i=1}^2
\frac{f_i(\tau_i,\xi_i)}
{\lb\lambda_i\rb^{b}}
\right\|_{L^2} \\
&\leq c \left\|
I^{\frac{\alpha}{2}}_{\ast}(v_1,v_2)
\right\|_{L^2}\leq c\prod\limits_{i=1}^2 \|v_i\|_{X_{0,0,b}}
=c\prod\limits_{i=1}^2 \|f_i\|_{L^2}
\end{align*}
since $b>1/2$.
For $J_{11,1}$ we use \eqref{eq:freq_bound} and \eqref{eq:sigma} again and get
$$\|J_{11,1}\|_{L^2}\leq c\left\|
\int_{\ast}\chi_{D_{11}\cap A_1}
\frac{|\xi|^{\frac{\alpha}{2}}}
{\min_{i=1,2}\lb\sigma_i\rb^{\omega}}\lb\lambda\rb^{b'}
f_1(\tau_1,\xi_1)\lb\lambda_1\rb^{\omega-b}f_2(\tau_2,\xi_2)\lb\lambda_2\rb^{-b}
\right\|_{L^2}
$$
We may assume that $|\lambda_1|\geq 2 |\lambda|$, because otherwise the same
argument as for $J_{11,0}$ applies.
If $\lb\sigma_1\rb\leq
\lb\sigma_2\rb$ we have $\lb\lambda_1\rb^{\omega} \leq
\min_{i=1,2}\lb\sigma_i\rb^{\omega}$.
If we suppose that $\lb\sigma_2\rb\leq \lb\sigma_1\rb$
we see
\begin{equation*}
|\lambda_1|= |\tau_1-\xi_1|\xi_1|^{\alpha}|
= |\tau-\tau_2-\xi|\xi|^{\alpha}+\xi|\xi|^{\alpha}-\xi_1|\xi_1|^{\alpha}| \leq |\lambda|+16|\sigma_2|
\end{equation*}
since we are in region $D_{11}$. This implies $\lb\lambda_1\rb\leq
c\lb\sigma_2\rb$ and we also have
$$\lb\lambda_1\rb^{\omega} \leq c
\min_{i=1,2}\lb\sigma_i\rb^{\omega}.$$
Therefore,
\begin{equation*}
\|J_{11,1}\|_{L^2} \leq c\left\|
\int_{\ast}\chi_{D_{11}\cap A_1}
|\xi|^{\frac{\alpha}{2}}
\lb\lambda\rb^{b'}
f_1(\tau_1,\xi_1)\lb\lambda_1\rb^{-b}f_2(\tau_2,\xi_2)\lb\lambda_2\rb^{-b}
\right\|_{L^2}
\end{equation*}
In $D_{11}$ we have $|\xi|^{\frac{\alpha}{2}}\leq c |
|\xi_2|^{\alpha}-|\xi_1|^{\alpha}|^{\frac{1}{2}}$ and by assumption
$b'\leq 0$, such that we may proceed as above with $J_0$ and use the
estimate \eqref{eq:bil_str} to conclude
\begin{equation*}
\|J_{11,1}\|_{L^2}\leq c
\left\|
\int_{\ast}
| |\xi_2|^{\alpha}-|\xi_1|^{\alpha}|^{\frac{1}{2}}
\prod\limits_{i=1}^2
\frac{f_i(\tau_i,\xi_i)}
{\lb\lambda_i\rb^{b}}
\right\|_{L^2}\leq
c\prod\limits_{i=1}^2 \|f_i\|_{L^2}
\end{equation*}

For $J_{11,2}$, we have by \eqref{eq:freq_bound} and \eqref{eq:sigma}
$$
\|J_{11,2}\|_{L^2}\leq c\left\|
\int_{\ast}\chi_{D_{11}\cap A_2}
|\xi|^{\frac{\alpha}{2}}\lb\lambda\rb^{b'}
f_1(\tau_1,\xi_1)\lb\lambda_1\rb^{-b}f_2(\tau_2,\xi_2)\lb\lambda_2\rb^{\omega-b}
\right\|_{L^2}
$$
In $D_{11} \cap A_2$ we have
$$
|\xi|^{\frac{\alpha}{2}}\leq c ||\xi|^{\alpha}-|\xi_1|^{\alpha}|^{\frac{1}{2}} \;\text{and
}\lb\lambda_2\rb^{\omega-b}\leq \lb\lambda\rb^{\omega-b}
$$
such that, because of $b'+\omega \leq 0$,
\begin{align*}
\|J_{11,2}\|_{L^2}&
\leq c \left\|K^{\frac{\alpha}{2}}_{\ast}(\overline{v}_1,\mathcal{F}^{-1}f_2) \right\|_{X_{0,0,-b}}\\
&\leq c\|v_1\|_{X_{0,0,b}}\|\mathcal{F}^{-1}f_2\|_{L^2}
=c\prod\limits_{i=1}^2 \|f_i\|_{L^2}
\end{align*}
for $b>1/2$ by the estimate \eqref{eq:dual_bil_str}.

Let us now consider the region $D_{12}$.
We define the contributions
$$
J_{12,0}=
\int_\ast\chi_{D_{12}\cap A}
|\xi|^{1-\omega}\lb\xi\rb^{s-\alpha\omega}\lb\lambda\rb^{b'}\lb\sigma\rb^{\omega}
\prod\limits_{i=1}^2
\frac{|\xi_i|^{\omega}\lb\xi_i\rb^{\alpha\omega-s}f_i(\tau_i,\xi_i)}
{\lb\lambda_i\rb^{b}\lb\sigma_i\rb^{\omega}}
$$
and, for $k=1,2$,
$$
J_{12,k}=
\int_\ast\chi_{D_{12}\cap A_k}
|\xi|^{1-\omega}\lb\xi\rb^{s-\alpha\omega}\lb\lambda\rb^{b'}\lb\sigma\rb^{\omega}
\prod\limits_{i=1}^2
\frac{|\xi_i|^{\omega}\lb\xi_i\rb^{\alpha\omega-s}f_i(\tau_i,\xi_i)}
{\lb\lambda_i\rb^{b}\lb\sigma_i\rb^{\omega}}
$$
In the subregion $D_{12}\cap A$ we use
$$ |\xi|^{-\alpha b'}\lb \xi_1 \rb^{-b'} \leq c
\lb\lambda\rb^{-b'}$$
and
\begin{align*}
&\|J_{12,0}\|_{L^2} \\& \leq c
\left\|\int_{\ast}\chi_{D_{12}\cap A}
|\xi|^{1-\omega+\alpha b'}\lb\xi\rb^{s-\alpha\omega}
\lb\sigma\rb^{\omega}\lb\xi_1\rb^{b'+\alpha\omega-s}\lb\xi_2\rb^{\alpha\omega-s}\prod\limits_{i=1}^2
\frac{f_i(\tau_i,\xi_i)|\xi_i|^{\omega}}
{\lb\lambda_i\rb^{b}\lb\sigma_i\rb^{\omega}}\right\|_{L^2}
\end{align*}
Using $\lb\xi_2\rb^{\alpha\omega-s}\lb\xi\rb^{s-\alpha\omega} \leq c$ and \eqref{eq:sigma} this is bounded by
$$
\left\|\int_{\ast}\chi_{D_{12}\cap A}
|\xi|^{1+\alpha b'}
\frac{\lb\xi_1\rb^{b'+\alpha\omega-s+\omega}}{\min_{i=1,2}\lb\sigma_i\rb^\omega}
\prod\limits_{i=1}^2\frac{f_i(\tau_i,\xi_i)}
{\lb\lambda_i\rb^{b}}\right\|_{L^2}
$$
Now, for $b'+\omega \leq 0$ we estimate
$$
|\xi|^{1+\alpha b'} \leq c
|\xi|^{\frac{\alpha}{2}}\lb\xi_1\rb^{1-\frac{\alpha}{2}+\alpha b'}
$$
since $1-\frac{\alpha}{2}+\alpha b'\leq 0$.
Moreover,
$$1-\frac{\alpha}{2}+\alpha b'+b'+\alpha\omega-s+\omega-(1+\alpha)\omega=
1-\frac{\alpha}{2}+\alpha b'+b'-s$$
which is negative for
$$s \geq \alpha (-\frac{1}{2}+b')+1+b'$$
Therefore, choosing $b'\leq \min\{-\omega,-\frac{1}{4}\}$, we continue for $s\geq
-\frac{3}{4}(\alpha-1)$ with
$$
\left\|\int_{\ast}\chi_{D_{12}\cap A}
|\xi|^{\frac{\alpha}{2}}
\prod\limits_{i=1}^2\frac{f_i(\tau_i,\xi_i)}
{\lb\lambda_i\rb^{b}}\right\|_{L^2}\leq c  \left\|
I^{\frac{\alpha}{2}}_{\ast}(v_1,v_2)
\right\|_{L^2}\leq c\prod\limits_{i=1}^2 \|f_i\|_{L^2}
$$
Next, we study the contribution of $J_{12,1}$. We may assume that
$\lb \lambda_1\rb \geq 2 \lb\lambda\rb$, because otherwise we use the same argument as in
$D_{12}\cap A$.
In $D_{12}\cap A_1$ we exploit
$$ |\xi|\lb \xi_1 \rb^{\frac{1}{\alpha}} \leq c
\lb\lambda_1\rb^{\frac{1}{\alpha}}$$
We observe that
$$|\lambda_1|=|\tau_1-\xi_1|\xi_1|^{\alpha}|\leq |\lambda|+c\lb\sigma_2\rb
\Rightarrow \lb \lambda_1\rb \leq c\lb\sigma_2\rb $$
and therefore
$$\lb \lambda_1 \rb^{\omega}\leq c\min_{i=1,2} \lb\sigma_i\rb^{\omega}$$
This shows
\begin{align*}
\|J_{12,1}\|_{L^2}& 
\leq c\left\|
\int_\ast\chi_{D_{12}\cap A_1}
\lb\lambda\rb^{b'}
\lb\xi_1\rb^{-\frac{1}{\alpha}+\omega+\alpha\omega-s}
\lb\lambda_1\rb^{\frac{1}{2}-b}\lb\lambda_2\rb^{-b}
\prod\limits_{i=1}^2
f_i(\tau_i,\xi_i)
\right\|_{L^2}
\end{align*}
We choose $b>\frac{1}{2}$ and in $D_{12}$ we have
$|\xi_1|\leq |\xi_2|$. Since we only consider $s\leq
\frac{1}{2}-\frac{\alpha}{2}$ (which means $\eps \leq \frac{\alpha-1}{4}$), we have
\begin{align*}
\|J_{12,1}\|_{L^2}& 
\leq c\left\|
\int_\ast
\lb\lambda\rb^{b'}
\lb\xi_2\rb^{\frac{1}{2}-\frac{\alpha}{2}-s}
\lb\lambda_2\rb^{-b}
\prod\limits_{i=1}^2
f_i(\tau_i,\xi_i)
\right\|_{L^2}
\end{align*}
With $b'\leq -\frac{1}{4}$ and Sobolev in time we see
\begin{align*}
\|J_{12,1}\|_{L^2}& 
\leq c\left\|\mathcal{F}^{-1}f_1
  J^{\frac{1}{2}-\frac{\alpha}{2}-s}v_2\right\|_{L^{4/3}_tL^2_x}\\
& \leq c \|f_1\|_{L^2_tL^2_x}
  \|J^{\frac{1}{2}-\frac{\alpha}{2}-s}v_2\|_{L^4_t L^{\infty}_x}
\end{align*}
Finally, by \eqref{eq:str_with_deriv}
$$
\|J^{\frac{1}{2}-\frac{\alpha}{2}-s}v_2\|_{L^4_t L^{\infty}_x}\leq c\|v_2\|_{X_{0,0,b}}=\|f_2\|_{L^2}
$$
if $\frac{1}{2}-\frac{\alpha}{2}-s \leq \frac{\alpha-1}{4}$, which is
equivalent to $s \geq -\frac{3}{4}(\alpha-1)$.

Now we turn to the contribution of $D_{12}\cap A_2$, where we use
$$ |\xi|^{-\alpha b'}\lb \xi_1 \rb^{-b'} \leq c
\lb\lambda_2\rb^{-b'}$$
and it follows
\begin{align*}
&\|J_{12,2}\|_{L^2} \\& \leq c
\left\|\int_{\ast}\chi_{D_{12}\cap A_2}
|\xi|^{1+\alpha
  b'}\frac{\lb\xi_1\rb^{b'+\omega+\alpha\omega-s}}{\min_{i=1,2}\lb\sigma_i\rb^{\omega}}
\lb\lambda\rb^{b'}\lb\lambda_2\rb^{-b'-b}\lb\lambda_1\rb^{-b}
\prod\limits_{i=1}^2 f_i(\tau_i,\xi_i)\right\|_{L^2}
\end{align*}
We have $$\lb\lambda\rb^{b'}\lb\lambda_2\rb^{-b'-b}\leq \lb\lambda\rb^{-b}$$ and
$$\min_{i=1,2}\lb\sigma_i\rb^{\omega}\geq \lb\xi_1\rb^{(1+\alpha)\omega}$$
and if $b'\leq -\omega$ we have $1+\alpha b'-\frac{\alpha}{2}\leq 0$ and therefore
$$|\xi|^{1+\alpha b'} \leq c |\xi|^{\frac{\alpha}{2}}\lb\xi_1\rb^{1+\alpha
  b'-\frac{\alpha}{2}}$$
If $b'\leq -\frac{1}{4}$ and $s\geq -\frac{3}{4}(\alpha-1)$
we estimate $b'-s+1+\alpha b'-\frac{\alpha}{2} \leq 0$ and
$$|\xi|^{\frac{\alpha}{2}}\leq c ||\xi|^{\alpha}-|\xi_1|^{\alpha}|^{\frac{1}{2}}$$
and therefore, by the dual bilinear Strichartz estimate \eqref{eq:dual_bil_str}
\begin{align*}
\|J_{12,2}\|_{L^2} & \leq c
\left\|\int_{\ast}
||\xi|^{\alpha}-|\xi_1|^{\alpha}|^{\frac{1}{2}}
\lb\lambda\rb^{-b}\lb\lambda_1\rb^{-b}
\prod\limits_{i=1}^2 f_i(\tau_i,\xi_i)\right\|_{L^2}\\
& \leq c \prod\limits_{i=1}^2 \|f_i\|_{L^2}
\end{align*}
This completes the discussion of the subregion $D_1$.

Let us now consider the
domain $D_2$, where $|\xi_1|\leq |\xi_2|\leq 4|\xi_1|$, $|\xi| \leq 2|\xi_2|$
and $|\xi|\leq 5|\xi_1|$. We subdivide $D_2=D_{21}\cup D_{22}$, where in
$$D_{21}\,:\; \xi_1\xi_2>0 \text{ or } |\xi| \geq \frac{1}{2}|\xi_1| \text{ or }
|\xi_2| \leq 1$$
and in
$$D_{22}\,:\; \xi_1\xi_2<0 \text{ and }|\xi| \leq \frac{1}{2}|\xi_1| \text{ and }
|\xi_2| \geq 1$$
additionally hold.
As above, we
define for $j=1,2$
$$
J_{2j,0}=
\int_\ast\chi_{D_{2j}\cap A}
|\xi|^{1-\omega}\lb\xi\rb^{s-\alpha\omega}\lb\lambda\rb^{b'}\lb\sigma\rb^{\omega}
\prod\limits_{i=1}^2
\frac{|\xi_i|^{\omega}\lb\xi_i\rb^{\alpha\omega-s}f_i(\tau_i,\xi_i)}
{\lb\lambda_i\rb^{b}\lb\sigma_i\rb^{\omega}}
$$
and for $k=1,2$
$$
J_{2j,k}=
\int_\ast\chi_{D_{2j}\cap A_k}
|\xi|^{1-\omega}\lb\xi\rb^{s-\alpha\omega}\lb\lambda\rb^{b'}\lb\sigma\rb^{\omega}
\prod\limits_{i=1}^2
\frac{|\xi_i|^{\omega}\lb\xi_i\rb^{\alpha\omega-s}f_i(\tau_i,\xi_i)}
{\lb\lambda_i\rb^{b}\lb\sigma_i\rb^{\omega}}
$$
We start with the
discussion of $D_{21}$, where all frequencies are of comparable size or smaller
then a constant, which shows that
$$
\frac{|\xi|^{1-\omega}\lb\xi\rb^{s-\alpha\omega}|\xi_1|^{\omega}|\xi_2|^{\omega}\lb\sigma\rb^{\omega}}
{\lb\xi_1\rb^{s-\alpha\omega}\lb\xi_2\rb^{s-\alpha\omega}\lb\sigma_1\rb^{\omega}\lb\sigma_2\rb^{\omega}}
\leq c \lb\xi\rb^{1-s}
$$
Therefore,
\begin{align*}
\|J_{21,0}\|_{L^2} & \leq c \left\|\int_\ast\chi_{D_{21}\cap A}
\lb\xi\rb^{1-s}\lb\lambda\rb^{b'}
\prod\limits_{i=1}^2
\frac{f_i(\tau_i,\xi_i)}
{\lb\lambda_i\rb^{b}}\right\|_{L^2}
\end{align*}
In $A$ we have $$\lb\xi\rb^{-b'(1+\alpha)}\leq c \lb\lambda\rb^{-b'}$$ and we use
the Strichartz estimate \eqref{eq:str_with_deriv} to conclude
\begin{align*}
\|J_{21,0}\|_{L^2} & \leq c \left\|\int_\ast
\lb\xi\rb^{1-s+b'(1+\alpha)-\frac{\alpha-1}{4}}
\lb\xi_1\rb^{\frac{\alpha-1}{4}}\mathcal{F}v_1(\tau_1,\xi_1)\mathcal{F}v_2(\tau_2,\xi_2)\right\|_{L^2}\\
 & \leq c\|J^{\frac{\alpha-1}{4}}v_1\|_{L_{t}^4L^{\infty}_x}\|v_2\|_{L^4_tL^{2}_x} \leq c\prod\limits_{i=1}^2 \|f_i\|_{L^2}
\end{align*}
since $1-s+b'(1+\alpha)-\frac{\alpha-1}{4}\leq 0$, which is equivalent
to $\frac{5}{4}+b'-\frac{\alpha}{4}+\alpha b'\leq s$. This is fulfilled for $b'
\leq -\frac{1}{2}+\frac{\eps}{3}$.
In $A_1$ we have $$\lb\xi\rb^{b(1+\alpha)}\leq c \lb\lambda_1\rb^{b}$$ and we
use Sobolev in time and the Strichartz estimate \eqref{eq:str_with_deriv} to
conclude for $b'\leq -\frac{1}{4}$
\begin{align*}
\|J_{21,1}\|_{L^2} & \leq c \left\|\int_\ast
\lb\xi\rb^{1-s-b(1+\alpha)-\frac{\alpha-1}{4}}\lb\lambda\rb^{b'}
f_1(\tau_1,\xi_1)\lb\xi_2\rb^{\frac{\alpha-1}{4}}\mathcal{F}v_2(\tau_2,\xi_2)\right\|_{L^2}\\
 & \leq c\|\mathcal{F}^{-1}f_1J^{\frac{\alpha-1}{4}}v_2\|_{L_{t}^{4/3}L^{2}_x}
\leq
 c\|f_1\|_{L_{tx}^2}\|J^{\frac{\alpha-1}{4}}v_2\|_{L^4_tL^{\infty}_x}\\
&\leq c\prod\limits_{i=1}^2 \|f_i\|_{L^2}
\end{align*}
The same argument applies to $J_{21,2}$ by exchanging the roles of $f_1,f_2$.

Finally, we turn to the contributions from the region $D_{22}$.
Here, we have $\xi_1\xi_2<0$. Therefore, we may write $\xi_1=\beta\xi_2$ for
$\beta \in [-1,-\frac{1}{4}]$. By the mean value theorem, this shows
\begin{equation}
  \label{eq:smoothing}
  \left||\xi_1|^{\alpha}-|\xi_2|^{\alpha}\right|^{\frac{1}{2}}
=\left||\beta|^{\alpha}-1\right|^{\frac{1}{2}}|\xi_2|^{\frac{\alpha}{2}}\geq
\frac{1}{2}||\beta|-1|^{\frac{1}{2}}|\xi_2|^{\frac{\alpha}{2}}=
\frac{1}{2}|\xi|^{\frac{1}{2}}|\xi_2|^{\frac{\alpha-1}{2}}
\end{equation}
Let us start with the subregion $A$.
We have $$\lb\sigma\rb^{\omega}\leq c \lb\lambda\rb^{\omega}+c\chi_{|\xi|\geq 1}\lb\xi\rb^{\omega+\alpha\omega}$$
which shows
\begin{align*}
\|J_{22,0}\|_{L^2}\leq & c
\left\|
\int_\ast\chi_{D_{22}\cap A}
|\xi|^{1-\omega}\lb\xi\rb^{s-\alpha\omega}\lb\lambda\rb^{b'+\omega}
\prod\limits_{i=1}^2
\frac{|\xi_i|^{\omega}\lb\xi_i\rb^{\alpha\omega-s}f_i(\tau_i,\xi_i)}
{\lb\lambda_i\rb^{b}\lb\sigma_i\rb^{\omega}}\right\|_{L^2}\\
& + c\left\|
\int_\ast\chi_{D_{22}\cap A}\chi_{|\xi| \geq 1}
\lb\xi\rb^{1+s}\lb\lambda\rb^{b'}
\prod\limits_{i=1}^2
\frac{|\xi_i|^{\omega}\lb\xi_i\rb^{\alpha\omega-s}f_i(\tau_i,\xi_i)}
{\lb\lambda_i\rb^{b}\lb\sigma_i\rb^{\omega}}\right\|_{L^2}
\end{align*}
Using
\begin{align*}
|\xi|^{-b'-\omega}\lb\xi_2\rb^{-\alpha b'-\alpha\omega} \leq c\lb\lambda\rb^{-b'-\omega}
\end{align*}
and \eqref{eq:smoothing}
we see that the first term is bounded by
\begin{align*}
&
\left\|
\int_\ast\chi_{D_{22}\cap A}
|\xi|^{1+b'}\lb\xi\rb^{s-\alpha\omega} \lb\xi_2\rb^{-2s+\alpha b'+\alpha \omega}
\prod\limits_{i=1}^2
\frac{f_i(\tau_i,\xi_i)}
{\lb\lambda_i\rb^{b}}\right\|_{L^2}\\
\leq &c \left\|
\int_\ast
\lb\xi\rb^{\frac{1}{2}+b'+s-\alpha\omega} \lb\xi_2\rb^{-2s+\alpha b'+\alpha \omega-\frac{\alpha-1}{2}}
\left||\xi_1|^{\alpha}-|\xi_2|^{\alpha}\right|^{\frac{1}{2}}
\prod\limits_{i=1}^2
\frac{f_i(\tau_i,\xi_i)}
{\lb\lambda_i\rb^{b}}\right\|_{L^2}
\end{align*}
If $b' \leq -\frac{1}{4}$ and $\eps \leq \frac{1}{4}$, then
$\frac{1}{2}+b'+s-\alpha\omega\leq 0$. Moreover, for $b' \leq
-\frac{1}{2}+\eps$, we have
$-2s+\alpha b'+\alpha \omega-\frac{\alpha-1}{2}\leq 0$. Then, by the bilinear
Strichartz estimate \eqref{eq:bil_str} this is bounded by
\begin{equation*}
\ldots \leq c\left\|
\int_\ast
\left||\xi_1|^{\alpha}-|\xi_2|^{\alpha}\right|^{\frac{1}{2}}
\prod\limits_{i=1}^2
\frac{f_i(\tau_i,\xi_i)}
{\lb\lambda_i\rb^{b}}\right\|_{L^2}
\leq c\prod\limits_{i=1}^2 \|f_i\|_{L^2}
\end{equation*}
For the second term we use
\begin{align*}
|\xi|^{-b'}\lb\xi_2\rb^{-\alpha b'} \leq c\lb\lambda\rb^{-b'}
\end{align*}
and find with \eqref{eq:smoothing}
\begin{align*}
\ldots &\leq c\left\|
\int_\ast\chi_{D_{22}\cap A} \chi_{|\xi| \geq 1}
\lb\xi\rb^{s+1+b'}\lb\xi_2\rb^{\alpha b'-2s}
\prod\limits_{i=1}^2
\frac{f_i(\tau_i,\xi_i)}
{\lb\lambda_i\rb^{b}}\right\|_{L^2}\\
& \leq c\left\|
\int_\ast 
\lb\xi\rb^{\frac{1}{2}+s+b'}\lb\xi_2\rb^{\alpha b'-2s-\frac{\alpha-1}{2}}
\left||\xi_1|^{\alpha}-|\xi_2|^{\alpha}\right|^{\frac{1}{2}}
\prod\limits_{i=1}^2
\frac{f_i(\tau_i,\xi_i)}
{\lb\lambda_i\rb^{b}}\right\|_{L^2}
\end{align*}
We only consider $\eps < \frac{3}{4}(\alpha-1)$. Then,
for $b' \leq -\frac{1}{2}+\frac{3}{4}(\alpha-1)-\eps$ we observe that
$\frac{1}{2}+s+b'\leq 0$.
Moreover $\alpha b'-2s-\frac{\alpha-1}{2}\leq 0$ for
$b' \leq -\frac{1}{2}+\eps$. Using the bilinear Strichartz estimate
\eqref{eq:bil_str}, we arrive at
$$\|J_{22,0}\|_{L^2}\leq c\prod\limits_{i=1}^2 \|f_i\|_{L^2}$$

Next, we consider the subregion $A_1$.
We have $$\lb\sigma\rb^{\omega}\leq c \lb\lambda_1\rb^{\omega}+c\chi_{|\xi|\geq 1}\lb\xi\rb^{\omega+\alpha\omega}$$
which shows
\begin{align*}
\|J_{22,1}\|_{L^2}\leq & c
\left\|
\int_\ast\chi_{D_{22}\cap A_1}
|\xi|^{1-\omega}\lb\xi\rb^{s-\alpha\omega}\lb\xi_2\rb^{-2s}
\lb\lambda\rb^{b'}\lb\lambda_1\rb^{-b+\omega}\lb\lambda_2\rb^{-b}
\prod\limits_{i=1}^2
f_i(\tau_i,\xi_i)\right\|_{L^2}\\
& + c\left\|
\int_\ast\chi_{D_{22}\cap A_1}\chi_{|\xi| \geq 1}
\lb\xi\rb^{1+s}\lb\lambda\rb^{b'}\lb\lambda_1\rb^{-b}\lb\lambda_2\rb^{-b}\lb\xi_2\rb^{-2s}
\prod\limits_{i=1}^2
f_i(\tau_i,\xi_i)\right\|_{L^2}
\end{align*}
As above, by
\begin{align*}
|\xi|^{-b'-\omega}\lb\xi_2\rb^{-\alpha b'-\alpha\omega} \leq c\lb\lambda_1\rb^{-b'-\omega}
\end{align*}
we see that the first term is bounded by
\begin{align*}
&
\left\|
\int_\ast\chi_{D_{22}\cap A}
|\xi|^{1+b'}\lb\xi\rb^{s-\alpha\omega} \lb\xi_2\rb^{-2s+\alpha b'+\alpha \omega}\lb\lambda\rb^{-b}\lb\lambda_2\rb^{-b}
\prod\limits_{i=1}^2
f_i(\tau_i,\xi_i)\right\|_{L^2}\\
& \leq c\left\|
\int_\ast
\lb\xi\rb^{1+b'+s-\alpha\omega} \lb\xi_2\rb^{-2s+\alpha b'+\alpha \omega-\frac{\alpha}{2}}
\left||\xi|^{\alpha}-|\xi_2|^{\alpha}\right|^{\frac{1}{2}}
\lb\lambda\rb^{-b}\lb\lambda_2\rb^{-b}
\prod\limits_{i=1}^2
f_i(\tau_i,\xi_i)\right\|_{L^2}
\end{align*}
Here, we used that due to $|\xi| \leq \frac{3}{4}|\xi_2|$ and $|\xi_2| \geq 1$ we have
$$\left||\xi|^{\alpha}-|\xi_2|^{\alpha}\right|^{\frac{1}{2}} \geq c\lb\xi_2\rb^{\frac{\alpha}{2}}$$
By estimating $\lb\xi\rb^{\frac{1}{2}} \leq \lb\xi_2\rb^{\frac{1}{2}}$ and with the
same restrictions on $s,b'$ as above we may apply the dual bilinear
Strichartz estimate \eqref{eq:dual_bil_str} and get
\begin{equation*}
\ldots \leq c\left\|
\int_\ast
\left||\xi|^{\alpha}-|\xi_2|^{\alpha}\right|^{\frac{1}{2}}
\lb\lambda\rb^{-b}\lb\lambda_2\rb^{-b}
\prod\limits_{i=1}^2
f_i(\tau_i,\xi_i)\right\|_{L^2}
\leq c\prod\limits_{i=1}^2 \|f_i\|_{L^2}
\end{equation*}
For the second term we use
\begin{align*}
|\xi|^{-b'}\lb\xi_2\rb^{-\alpha b'} \leq c\lb\lambda_1\rb^{-b'}
\end{align*}
and find
\begin{align*}
\ldots &\leq c\left\|
\int_\ast 
\lb\xi\rb^{1+s+b'}\lb\xi_2\rb^{\alpha b'-2s-\frac{\alpha}{2}}
\left||\xi|^{\alpha}-|\xi_2|^{\alpha}\right|^{\frac{1}{2}}
\lb\lambda\rb^{-b}\lb\lambda_2\rb^{-b}
\prod\limits_{i=1}^2
f_i(\tau_i,\xi_i)
\right\|_{L^2}\\
&\leq c\prod\limits_{i=1}^2 \|f_i\|_{L^2}
\end{align*}
by \eqref{eq:dual_bil_str} with the same restrictions on $s,b',b$ as in the
region $A$, since $\lb\xi\rb^{\frac{1}{2}} \leq \lb\xi_2\rb^{\frac{1}{2}}$.

Finally, we turn to the region $A_2$. In $D_{22}$ 
the frequencies $\xi_1$ and $\xi_2$ are of comparable size and due to
$|\xi| \leq \frac{1}{2}|\xi_1|$ and $|\xi_1| \geq
\frac{1}{4}|\xi_2|\geq\frac{1}{4}$ we have
$$\left||\xi|^{\alpha}-|\xi_1|^{\alpha}\right|^{\frac{1}{2}} \geq c\lb \xi_1\rb^{\frac{\alpha}{2}}$$
Now we use the same argument as $A_1$ with the roles of $f_1,f_2$ exchanged.

This finishes
the proof of the bilinear estimate for $s=s_0=-\frac{3}{4}(\alpha-1)+\eps$, for
$\eps \leq \frac{\alpha-1}{4}$. The restrictions on $b'$ can be summarized to
$$b' \leq
\min\{-\frac{1}{4},-\omega,-\frac{1}{2}+\frac{\eps}{3},-\frac{1}{2}+\frac{3}{4}(\alpha-1)-\eps\}$$
For $b$ we assumed $\frac{1}{2}<b<b'+1$.
Now we turn to the case $s >s_0=-\frac{3}{4}(\alpha-1)+\eps$. Let $\rho=s
-s_0$. Because of
$$\lb\xi\rb^\rho \leq c\lb\xi_1\rb^\rho+c\lb\xi_2\rb^\rho$$
we see
\begin{align*}
\|\partial_x(u_1u_2)\|_{X_{s,\omega,b'}}
&\leq c \|\partial_x(J^{\rho}u_1u_2)\|_{X_{s_0,\omega,b'}}+
\|\partial_x(u_1J^{\rho}u_2)\|_{X_{s_0,\omega,b'}}\\
&\leq c \|u_1\|_{X_{s,\omega,b}}\|u_2\|_{X_{s_0,\omega,b}}
+\|u_1\|_{X_{s_0,\omega,b}}\|u_2\|_{X_{s,\omega,b}}
\end{align*}
This proves that for all $s \geq s_0>-\frac{3}{4}(\alpha-1)$ we find suitable
$b' \in (-\frac{1}{2},0)$ and $b \in (\frac{1}{2},b'+1)$ such that the bilinear
estimate holds true. \hfill \qed

\section{An a priori bound}\label{sect:apr}
This section is devoted to the proof of an a priori bound for the
$H^{(0,\omega)}$ norm, which allows an iteration of the local argument to
prove global well-posedness for $s\geq 0$.
\begin{lemma}\label{lem:apriori_bound}
Let $s\geq 0$. There exists $C>0$, such that for all
smooth, real valued solutions $u$ of \eqref{eq:bo_alpha}, we have
\begin{equation}\label{eq:apriori_bound}
\sup_{t\in [-T,T]}\|u(t)\|_{H^{(0,\omega)}}
\leq C\|u(0)\|_{H^{(0,\omega)}}+CT\|u(0)\|^2_{H^{(0,\omega)}}
\end{equation}
\end{lemma}
\begin{proof}
We easily verify the conservation law
$$\|u(t)\|^2_{L^2}=\|u(0)\|^2_{L^2}, \quad t \in (-T,T)$$
Therefore it suffices to prove an a priori estimate for the low frequency part in $\dot{H}^{-\omega}$.
Let $\psi \in C^{\infty}_{0}([-2,2])$ be nonnegative with $\psi|_{[-1,1]} \equiv 1$.
We define $$\mathcal{F}_x v(t)(\xi)=\psi(\xi)|\xi|^{-\omega}\mathcal{F}_x u(t)(\xi)$$
The function $v$ solves the equation
\begin{align*}
v_t-|D|^{\alpha}v_x&=f \;\text{ in } (-T,T) \times \R\\
v(0)&=v_0
\end{align*}
where $\mathcal{F}_x v_0(\xi)=\psi(\xi)|\xi|^{-\omega}\mathcal{F}_x u(0)(\xi)$ and
$$\mathcal{F}_x f(t)(\xi)=-\frac{i}{2}\psi(\xi)\xi|\xi|^{-\omega}\mathcal{F}_x u^2(t)(\xi)$$
For fixed $t$ we estimate
\begin{align*}
\|f(t)\|_{L^2_x} \leq&
c\|\psi(\xi)\mathcal{F}_x u^2(t)(\xi)\|_{L^2_{\xi}}\leq
c\|\mathcal{F}_x u^2(t)\|_{L^{\infty}_{\xi}}\\
\leq & c\|u^2(t)\|_{L^1_{x}} \leq c\|u(t)\|^2_{L^2_{x}}
\end{align*}
This shows
\begin{align*}
\|v\|_{L^{\infty}_TL^2_x} &\leq c \|v_0\|_{L^2_x}+c\|f\|_{L^1_TL^2_x}
\leq c\|u(0)\|_{H^{(0,\omega)}}+cT\|u\|^2_{L^{\infty}_TL^2_x}\\
&\leq c\|u(0)\|_{H^{(0,\omega)}}+cT\|u(0)\|^2_{H^{(0,\omega)}}
\end{align*}
\end{proof}

\bibliographystyle{hplain}
\bibliography{literatur}\label{sect:refs}

\begin{thebibliography}{10}

\bibitem{Bo93}
J.~Bourgain.
\newblock Fourier transform restriction phenomena for certain lattice subsets
  and applications to nonlinear evolution equations.
\newblock {\em Geom. Funct. Anal.}, 3(2/3):107--156/209--262, 1993.

\bibitem{BP05}
N.~Burq and F.~Planchon.
\newblock {On well-posedness for the Benjamin-Ono equation}, 2005,
  arXiv:math.AP/0509096.

\bibitem{CKS02}
J.~Colliander, C.~Kenig, and G.~Staffilani.
\newblock Local well-posedness for dispersion-generalized {B}enjamin-{O}no
  equations.
\newblock {\em Differential Integral Equations}, 16(12):1441--1472, 2003.

\bibitem{G96}
J.~Ginibre.
\newblock {Le probl{\`e}me de Cauchy pour des EDP semi-lin{\'e}aires
  p{\'e}riodiques en variables d'espace [d'apr{\`e}s Bourgain].}
\newblock In {\em {S{\'e}minaire Bourbaki. Volume 1994/95. Expos{\'e}s
  790-804}}, pages {163--187}. {Soci{\'e}t{\'e} Math{\'e}matique de France.
  Paris: Ast{\'e}risque. 237, Exp. No.796}, 1996.

\bibitem{G01b}
A.~Gr{\"u}nrock.
\newblock {A bilinear Airy- estimate with application to gKdV-3}, 2001,
  arXiv:math.AP/0108184.

\bibitem{G01a}
A.~Gr{\"u}nrock.
\newblock {Some local wellposedness results for nonlinear Schr\"odinger
  equations below $L^2$}, 2001, arXiv:math.AP/0011157.

\bibitem{H05a}
S.~Herr.
\newblock {Well-posedness for equations of Benjamin-Ono type}.
\newblock Submitted for publication, 2005,
  http://www.mathematik.uni-dortmund.de/$\sim$herr.

\bibitem{IK05}
A.~Ionescu and C.E. Kenig.
\newblock {Global well-posedness of the Benjamin-Ono equation in low-regularity
  spaces}, 2005, arXiv:math.AP/0508632.

\bibitem{KK03}
C.E. Kenig and K.D. Koenig.
\newblock {O}n the local well-posedness of the {B}enjamin-{O}no and modified
  {B}enjamin-{O}no equations.
\newblock {\em Math. Res. Lett.}, 10(5-6):879--895, 2003.

\bibitem{KPV91b}
C.E. Kenig, G.~Ponce, and L.~Vega.
\newblock Oscillatory integrals and regularity of dispersive equations.
\newblock {\em Indiana Univ. Math. J.}, 40(1):33--69, 1991.

\bibitem{KPV96}
C.E. Kenig, G.~Ponce, and L.~Vega.
\newblock A bilinear estimate with applications to the {K}d{V} equation.
\newblock {\em J. Amer. Math. Soc.}, 9(2):573--603, 1996.

\bibitem{MR03b}
L.~Molinet and F.~Ribaud.
\newblock On global well-posedness for a class of nonlocal dispersive wave
  equations.
\newblock Preprint.

\bibitem{MST01}
L.~Molinet, J.C. Saut, and N.~Tzvetkov.
\newblock Ill-posedness issues for the {B}enjamin-{O}no and related equations.
\newblock {\em SIAM J. Math. Anal.}, 33(4):982--988 (electronic), 2001.

\bibitem{T03}
T.~Tao.
\newblock {Global well-posedness of the Benjamin-Ono equation in $H^1(\mathbb
  R)$.}
\newblock {\em J. Hyperbolic Differential Equations}, 1:27--49, 2004,
  arXiv:math.AP/0307289.

\end{thebibliography}

\end{document}